\magnification=1200

\def\hal{{\vrule height 10pt width 4pt depth 0pt}}

\def\a{{\alpha}}
\def\b{{\beta}}
\def\B{{\cal B}}
\def\C{{\bf C}}
\def\f{{1\over 4}}
\def\ke{{\rm ker}}
\def\la{{\langle}}
\def\ra{{\rangle}}
\def\R{{\bf R}}
\def\Prim{{\rm Prim}}


\def\ep{{\varepsilon}}

\def\ld{{\lambda}}

\def\ph{{\varphi}}
\def\ps{{\psi}}
\def\rh{{\rho}}


\def\vphi{{\ph}}
\def\vepsilon{{\ep}}

\centerline{\bf Modules with norms which take values in a C*-algebra}
\medskip

\centerline{N.\ C.\ Phillips\footnote{*}{Partially
supported by NSF grant DMS-9400904.}
and N. Weaver\footnote{**}{Partially
supported by an NSF postdoctoral fellowship.}}
\bigskip
\bigskip

{\narrower{
\it We consider modules $E$ over a C*-algebra $A$ which are
equipped with a
map into $A_+$ that has the formal properties of a norm. We completely
determine the structure of these modules. In particular, we show
that if $A$ has no nonzero commutative ideals then
every such $E$ must be a
Hilbert module.
The commutative case is much less rigid:
if $A = C_0(X)$ is commutative then
$E$ is merely isomorphic to the module of continuous
sections of some bundle
of Banach spaces over $X$. In general $E$ will embed in a direct sum of
modules of the preceding two types.
\bigskip}}
\bigskip

Let $A$ be a C*-algebra, and let $A_+$ denote the set of positive
elements of $A.$ We define a
{\bf Finsler $A$-module} to be a left
$A$-module $E$ which is equipped with a map $\rho: E \to A_+$ such that
\medskip

{\narrower{
\noindent (1) the map $\|\cdot\|_E: x \mapsto \|\rho(x)\|$ is a Banach
space norm on $E$; and
\medskip

\noindent (2) $\rho(ax)^2 = a \rho(x)^2 a^*$ for all
$a \in A$ and $x \in E$.
\medskip}}

\noindent If we use the convention
$|b| = (b b^*)^{1/2}$ for $b \in A$, then
condition (2) is equivalent to
$$
\rho(ax) = |a \rho(x)|.
$$
For $A$ commutative this is the same as
$\rho(ax) = |a| \rho(x)$, which is the
usual form this sort of axiom takes in the commutative case.
But this last
version is not appropriate in the noncommutative case because
$\rho(ax)$ is
positive, while $|a| \rho(x)$, though a product of two positive
elements,
in general cannot be expected to be self-adjoint, let alone positive.
\medskip

If $E$ is a Hilbert $A$-module then defining
$\rho(x) = \la x, x\ra^{1/2}$
makes $E$ a Finsler module; in particular,
$\la ax, ax\ra = a \la x, x\ra a^*$, so condition (2) holds. This also
helps to justify the specific form of this condition,
on the grounds that any
definition of an $A$-valued norm ought to include norms arising from
Hilbert modules in this manner.
\medskip

Indeed, in the commutative case Finsler modules
are a natural generalization
of Hilbert modules. To see this let $X$ be a locally
compact space and let $\B = \bigcup_{t \in X} H_t$
be a bundle of Hilbert
spaces over $X$ satisfying appropriate continuity properties.
Then the set
$E$ of continuous sections (that is, continuous maps $f: X \to \B$ such
that $f(t) \in H_t$
for $t \in X$) which vanish at infinity, is naturally a $C_0(X)$-module.
Furthermore it has a $C_0(X)$-valued inner product defined by
$$
\la f, g\ra(t) = \la f(t), g(t)\ra_{H_t}
$$
for $t \in X$, hence is a Hilbert $C_0(X)$-module [16].
Conversely, every
Hilbert $C_0(X)$-module is isomorphic to one of this form [21].
\medskip

If we ask instead only that $\B = \bigcup_{t \in X} B_t$
be a bundle of Banach
spaces over $X$, then the module of continuous sections now possesses a
$C_0(X)$-valued norm
$$
\rho(f)(t) = \|f(t)\|_{B_t}
$$
rather than a $C_0(X)$-valued inner product.
It is easy to see that this makes
$E$ a Finsler $C_0(X)$-module, and we prove conversely that (as an easy
consequence of known facts) every Finsler
$C_0(X)$-module is isomorphic to one of this form.
\medskip

Thus, given the well-known conception of finitely generated projective
modules over C*-algebras as
``noncommutative vector bundles'' and Hilbert modules as
``noncommutative
Hilbert bundles'' ([18], [20]), it may appear that
our Finsler modules might serve as the basis for a noncommutative
version
of Banach bundles. Now we mentioned above that every Hilbert $A$-module
carries a natural Finsler structure.
One might hope to construct non-Hilbert Finsler modules over many
C*-algebras $A$ by forming a suitable completion of the algebraic
tensor product of $A$ with a non-Hilbert Banach space.
Surprisingly, we found that for ``most'' noncommutative C*-algebras,
namely all those algebras $A$ with no
nonzero commutative ideals, every Finsler $A$-module must arise
from a unique Hilbert $A$-module (Corollary 18).
In comparison with the commutative situation just discussed,
even with the case $A = \C$ (when $E$ can be any Banach space), the
noncommutative case is evidently far more rigid. From one standpoint
this is merely a negative result which shows that $A$-valued norms
are not interesting in the noncommutative case. On the other hand it
may be viewed as a positive result about the robustness of the concept
of Hilbert modules, a topic also explored in [8], and also as indirect
evidence that operator modules are really the right noncommutative
version of Banach bundles, a position we argue in section 2.
The Banach module properties of Hilbert modules have also been
considered in [14].
\medskip

Our terminology was chosen for the following reason.
A natural example of a bundle of Hilbert spaces is given by
the tangent bundle of a Riemannian manifold $X$.
Here the vector space over a point $t \in X$ is
simply the tangent space at $t$, and the fact that $X$ is Riemannian
means precisely that each tangent space has an inner product.
Finsler geometry is an increasingly popular generalization of
Riemannian geometry in which one requires only that each tangent
space have a norm ([5], [9]).
Thus Finsler geometry appears to involve Banach bundles
in the same way that Riemannian geometry involves Hilbert bundles.
We wish to thank David Blecher for pointing out this connection between
Finsler manifolds and Banach bundles.
\medskip

Section 1 contains preliminary general results.
In section 2 we establish
connections between operator modules, Finsler modules, and
Banach bundles in the commutative case. In section 3 we consider the
noncommutative case and obtain a complete description of the structure
of an arbitrary Finsler module.
\medskip

It is a pleasure to thank Charles Akemann for supplying a general
C*-algebra fact, Theorem 4.
This is a crucial result for our purposes and is also of
independent interest.
\bigskip
\bigskip

\noindent {\underbar{\bf 1. Preliminaries.}}
\bigskip

In this section we collect some important general facts about Finsler
modules. Aside from Akemann's result (Theorem 4) the material is fairly
trivial.
\medskip

Recall that a Banach $A$-module is an $A$-module $E$ that is
simultaneously a Banach space and which satisfies
$\|ax\| \leq \|a\| \|x\|$ for all $a \in A$ and $x \in E$.
\bigskip

\noindent {\bf Proposition 1.} {\it Every Finsler $A$-module is a Banach
$A$-module.}
\medskip

\noindent {\it Proof.}
By definition, $\|\cdot\|_E$ makes $E$ a Banach space.
We must therefore show that $\|ax\|_E \leq \|a\| \|x\|_E$ for
all $a \in A$ and $x \in E$. This follows from condition (2):
$$
\|ax\|_E^2 = \|\rho(ax)^2\| = \|a \rho(x)^2 a^*\|
\leq \|a\|^2 \|x\|_E^2.\eqno{\hal}
$$

In the next result we observe that if $A$ is commutative, every Finsler
$A$-module has properties which make it look very much like a
module with an $A$-valued norm.
Of these, the first has a natural analog in Finsler condition (2), as
we observed in the introduction.
\medskip

The second property, a generalized triangle inequality, is
far too strong in the noncommutative case (but see [2]),
although we see in this proposition
that if $A$ is commutative it follows from the seemingly weaker
assumption that $\|\cdot\|_E$ satisfies the triangle inequality.
The latter is suitable in the noncommutative setting, and
is already sufficient for the rather strong
structure results to be given in section 3.
On the other hand, in Lemma 12 we give a kind of noncommutative
generalization of this part of the proposition.
\medskip

\noindent {\bf Proposition 2.}
{\it Let $A$ be a C*-algebra with center $Z (A),$ and let
$E$ be a Finsler $A$-module. Then $\rho$ satisfies
$$
\rho(ax) = |a| \rho(x)
$$
for all $a \in Z(A)$ and $x \in E$. If $A = C_0(X)$ is commutative then
$\rho$ satisfies
$$
\rho(x + y) \leq \rho(x) + \rho(y)
$$
for all $x,y \in E$.}
\medskip

\noindent {\it Proof.}
If $a$ belongs to the center of $A$ then so does $|a|$,
hence both commute with $\rho(x)^2$ and
$$
\rho(ax)^2 = a \rho(x)^2 a^* = |a|^2 \rho(x)^2 = |a| \rho(x)^2 |a|;
$$
taking square roots yields $\rho(ax) = |a|\rho(x)$. To prove the second
statement suppose it is not the case and find a point $t \in X$
such that
$$
\rho(x + y)(t) > \a > \rho(x)(t) + \rho(y)(t).
$$
Let $U$ be a compact neighborhood of $t$ such that $\rho(x)(t') +
\rho(y)(t') < \a$ for all $t' \in U$.
Let $f \in C_0(X)$ satisfy $0 \leq f
\leq 1$, $f|_{X - U} = 0$, and $f(t) = 1$.
Then using the first part of this proposition we have
$$
\eqalign{\|fx\|_E + \|fy\|_E &= \|\rho(fx)\| + \|\rho(fy)\|
= \|f \rho(x)\| + \|f \rho(y)\| \leq \a\cr
& < f(t) \rho(x + y)(t) = \rho(fx + fy)(t)\cr
& \leq \|\rho(fx + fy)\| = \|fx + fy\|_E,\cr}
$$
contradicting the triangle inequality in $E$ (Finsler condition (1)).
This establishes the result.\hfill\hal
\bigskip

Next we observe that every Hilbert module gives rise to
a Finsler module.
\bigskip

\noindent {\bf Proposition 3.} {\it Let $A$ be a C*-algebra and $E$ a
Hilbert $A$-module. Then defining $\rho(x) = \la x, x\ra^{1/2}$
makes $E$ a Finsler $A$-module.}
\medskip

\noindent {\it Proof.} The fact that $\|\cdot\|_E$ is a complete norm is
part of the definition of a Hilbert module. The second Finsler condition
holds because
$$
\rho(ax)^2 =
   \la ax, ax\ra = a \la x, x\ra a^* = a \rho(x)^2 a^*.\eqno{\hal}
$$

Finally, we come to the one substantive result in this section [1].
It is
needed to show that $\rho$ is continuous and uniquely determined by the
scalar norm $\|\cdot\|_E$ (Corollaries 5 and 6), and also to show that
Finsler modules can be factored (Lemma 12).
\bigskip

\noindent {\bf Theorem 4 (Akemann).}
{\it Let $b$ and $c$ be positive elements
of a C*-algebra $A$. Then
$$
\|b - c\| = \sup\{\big| \|a b a\| - \|a c a\| \big|:
  a \in A_+, \, \|a\| \leq 1\}.
$$}

\noindent {\it Proof.} The inequality $\geq$ is easy since
$$
\big| \|a b a\| - \|a c a\| \big| \leq
       \|a b a - a c a\| \leq \|a\| \|b - c\| \|a\|
\leq \|b - c\|
$$
for any $a \in A_+$ with $\|a\| \leq 1$.
\medskip

Now for the reverse inequality. Without loss of generality assume that
$0 \leq b,c \leq 1$. By replacing $A$ with $C^*(b, c)$ we can
assume that $A$ is separable.  Assume that $\a = \|b - c\| > 0$.
It now suffices to prove that
there is a sequence $\{a_n\}$ in the positive unit ball of  
$A$ such that
$$
\lim_{n\to\infty} \big|\|a_n b a_n\| - \|a_n c a_n\|\big| = \a.
$$

To prove this, observe that,
exchanging $b$ and $c$ if necessary, there is a pure state $f$
of $A$ such that $f(b - c) = \a$. By Proposition 2.2 of [3], there is a
sequence $\{a_n\}$ of positive norm 1 elements in $A$ that excises $f$.
This means that for any $a$ in $A$, $\lim \|a_n a a_n - f(a) a_n^2\| =
0$.
Taking $a = b$ and then $a = c$, we get
$$
0 = \lim \|a_n b a_n - f(b) a_n^2\| \geq
\lim \big| \|a_n b a_n\| - \|f(b) a_n^2\| \big| \geq 0
$$
and
$$
0 = \lim \|a_n c a_n - f(c) a_n^2\| \geq
\lim \big| \|a_n c a_n\| - \|f(c) a_n^2\| \big| \geq 0.
$$
Since $\|f(b) a_n^2\|  = f(b)$ and $\|f(c) a_n^2\|  = f(c)$, this means
that $\lim \|a_n b a_n\| = f(b)$ and $\lim \|a_n c a_n\| = f(c)$.
Since $f(b) - f(c) = \a$, we get
$\lim (\|a_n b a_n\| - \|a_n c a_n\|) = \a,$ as desired.\hfill\hal
\bigskip

\noindent {\bf Corollary 5.} {\it Let $E$ be a Finsler module over a
C*-algebra $A$.
Then $\rho: E \to A_+$ is continuous from $\|\cdot\|_E$ to
$\|\cdot\|$.}
\medskip

\noindent {\it Proof.} Let $x_n, x \in E$ and set $b_n = \rho(x_n)^2$,
$b = \rho(x)^2$.
Suppose $\|x_n - x\|_E \to 0$ and let $C = \sup\{\|x_n\|_E\}$.
Then for any $a$ in the positive unit ball of $A$ we have
$$
\eqalign{\big| \|a b_n a\| - \|a b a\| \big|
&= \big| \|\rho(a x_n)^2\| - \|\rho(a x)^2\| \big|\cr
&= \big| \|a x_n\|_E^2 - \|a x\|_E^2 \big|\cr
&= \bigl( \|a x_n\|_E + \|a x\|_E \bigr)
            \big| \|a x_n\|_E - \|a x\|_E \big|\cr
&\leq 2C \|a(x_n - x)\|_E\cr
&\leq 2C \|x_n - x\|_E.\cr}
$$
Thus Theorem 4 implies
$$
\|\rho(x_n)^2 - \rho(x)^2\| = \|b_n - b\| \leq 2C \|x_n - x\|_E \to 0,
$$
so that we have $\rho(x_n)^2 \to \rho(x)^2$.
Finally, taking square roots implies that $\rho(x_n) \to \rho(x)$.
(Application of any continuous function preserves norm convergence of
self-adjoint elements in any C*-algebra; this
can be seen by first checking on polynomials.) Thus, $x_n \to x$ in $E$
implies $\rho(x_n) \to \rho(x)$ in $A$.\hfill\hal
\bigskip

\noindent {\bf Corollary 6.} {\it Let $A$ be a
C*-algebra and $E$ an $A$-module. Suppose $\rho, \rho': E \to A_+$
are two functions both of which make $E$
a Finsler module and which induce
the same norm $\|\cdot\|_E$. Then $\rho = \rho'$.}
\medskip

\noindent {\it Proof.}
Suppose $\rho(x)^2 \neq \rho'(x)^2$ for some $x \in E$. Then
$\|\rho(x)^2 - \rho'(x)^2\| \neq 0$ and so Theorem 4 implies that there
exists $a \in A_+$ with $\|a\| \leq 1$
such that $\|a \rho(x)^2 a\| \neq \|a \rho'(x)^2 a\|$.
Hence
$$
\|ax\|_E^2 = \|\rho(ax)^2\| \neq \|\rho'(ax)^2\| = \|ax\|_E^2,
$$
a contradiction.\hfill\hal
\bigskip

A similar statement about Hilbert modules, to the effect that the
scalar norm determines the
$A$-valued inner product, was given in [8].
This also follows from Corollary 6 since the $A$-valued inner
product is uniquely determined via polarization
from the Finsler norm that it gives rise to via Proposition 3.
\bigskip
\bigskip

\noindent {\underbar{\bf 2. Banach bundles.}}
\bigskip

A need for a noncommutative version of Banach bundles arises in the
theory of noncommutative metrics. This happens in the following way.
First of all, if $X$ is a Riemannian manifold then the cotangent bundle
is a Hilbert bundle, as mentioned in the introduction.
For $X$ any metric space there is a corresponding construction [11]
which involves Banach bundles, and this may be regarded as an
integrated version of the cotangent bundle construction [22].
This construction of de Leeuw actually ``encodes'' the metric structure
of $X$ in a manner so robust as to suggest that a notion of a
noncommutative
metric could be based on a noncommutative version of the set-up [22]. To
describe this noncommutative scheme one needs a noncommutative version
of the notion of a Banach bundle.
\medskip

On the basis of examples it has become clear that in describing
noncommutative metrics, Hilbert modules are sufficient for situations
in which
one has ``noncommutative Riemannian structure'' ([19], [23]), but more
generally one needs operator modules ([24], [25]).
Thus our first goal here is to show how in the
commutative case operator modules correspond to Banach
bundles, which suggests that general operator modules may be viewed as
noncommutative Banach bundles.
Modules associated to Banach bundles have been
thoroughly studied and so our results in this section are fairly easy
consequences of known facts.
\medskip

Before proceeding we must introduce a distinction emphasized in [12],
between {\bf (F) Banach bundles} and {\bf (H) Banach bundles}.
These are the bundle notions which respectively correspond to the
concepts of continuous fields of Banach spaces [15] and uniform fields
of Banach spaces [10]. In brief, the topology interacts with the norm
in such a way that the fiberwise norm of a continuous section of an
(F) Banach bundle is continuous, while in an (H) Banach bundle it need
only be semicontinuous.
\medskip

We also need the following definitions.
Let $A = C(X)$ be a unital commutative C*-algebra.
By an {\bf abelian operator $A$-module} (see [13]) we mean a Banach
$A$-module $E$ for which there exists a commutative C*-algebra $B$
together with an isometric embedding $\pi: E \to B$ and a $*$-isomorphic
embedding $\vphi: A \to B$, such that
$$
\pi(ax) = \vphi(a) \pi(x)
$$
for $a \in A$ and $x \in E$.
An {\bf $A$-convex $A$-module} [12] is a Banach
$A$-module $E$ which satsifies
$$
\|fx + gy\| \leq {\rm max}(\|x\|, \|y\|)
$$
for any $x,y \in E$ and any positive $f,g \in C(X)$ such that
$f + g = 1$.
\bigskip

\noindent {\bf Theorem 7.} {\it Let $A = C(X)$ be a unital commutative
C*-algebra and let $E$ be a Banach $A$-module.
The following are equivalent:
\medskip

{\narrower{
\noindent (a) $E$ is an abelian operator $A$-module.
\medskip

\noindent (b) $E$ is an $A$-convex $A$-module.
\medskip

\noindent (c) There is an (H) Banach bundle over $X$ of which
$E$ is
isomorphic to the module of continuous sections.
\medskip}}}

\noindent {\it Proof.} (a) $\Rightarrow$ (b). Let $B = C_0(Y)$ be a
commutative C*-algebra and suppose $A$ and $E$ are embedded in $B$.
Then the $A$-convex inequality is trivially checked at each $t \in Y$:
$$
|f(t) x(t) + g(t) y(t)| \leq (f(t) + g(t)) {\rm max}(\|x\|, \|y\|)
\leq {\rm max}(\|x\|, \|y\|),
$$
hence $\|fx + gy\| \leq {\rm max}(\|x\|, \|y\|)$.
\medskip

(b) $\Rightarrow$ (c). This is ([12], Theorem 2.5).
\medskip

(c) $\Rightarrow$ (a). Let $\B = \bigcup_{t \in X} B_t$ be an (H) Banach
bundle over $X$ and let $E$ be the module of continuous sections of
$\B$. Let
$$
Y = \{(t, v): t \in X, v \in B_t^*, \|v\| \leq 1\}
$$
where $B_t^*$ is the dual Banach space to $B_t$.
Then $A = C(X)$ embeds in $B = l^\infty(Y)$ by setting
$\vphi(f)(t, v) = f(t)$, and $E$ embeds in $B$
by setting $\pi(x)(t, v) = v(x(t))$. The module structure is preserved
by
these embeddings, for
$$
\pi(f x)(t, v) = v(f x(t)) = v(f(t) x(t)) = f(t) v(x(t)) =
                \vphi(f) \pi(x)(t, v).
$$
Thus $E$ is an abelian operator $A$-module.\hfill\hal
\bigskip

We view Theorem 7 as justifying the idea that general operator
modules are
``noncommutative Banach bundles.'' Note that the equivalence of parts
(a) and (c) easily extends to the case where $A = C_0(X)$ is nonunital,
since any Banach module over $A$ is also a Banach module over the
unitization of
$A$. But now part (c) will involve a Banach bundle over the one-point
compactification of $X$.
\medskip

Next we prove a similar fact about Finsler $C_0(X)$-modules.
\bigskip

\noindent {\bf Theorem 8.}
{\it Let $A = C_0(X)$ be a commutative C*-algebra
and let $E$ be a Banach $A$-module. The following are equivalent:
\medskip

{\narrower{
\noindent (a)
There exists a map $\rho: E \to A_+$ (necessarily unique) which induces
the given norm on $E$ and makes $E$ into a Finsler $A$-module.
\medskip

\noindent (b)
There is an (F) Banach bundle over $X$ of which $E$ is
isomorphic to the module of continuous sections vanishing at infinity.
\medskip}}}

\noindent {\it Proof.}
(a) $\Rightarrow$ (b). For $X$ compact (b) follows
from page 48 of [12] and Proposition 2.
If $X$ is locally compact but not compact,
let $X^+$ be its one-point
compactification; then $E$ is also a Finsler module over
$C(X^+)$. Hence $E$ is isomorphic
to the module of continuous sections of some (F) Banach bundle over
$X^+$.
But since $\rho(x) \in C_0(X)$ for all $x \in E$, the fiber over
$\infty$ must be trivial. So $E$ is also isomorphic to the module of
continuous sections vanishing at infinity of the restriction of the
bundle to $X$.
\medskip

(b) $\Rightarrow$ (a). We have already described the construction of
$\rho$ in the introduction, namely $\rho(x)(t) = \|x(t)\|_{B_t}$
for any section $x: X \to \B$.
This is a continuous function of $t$ precisely by the definition of an
(F) Banach bundle, and it satisfies the Finsler conditions because it
satisfies them fiberwise.
Uniqueness of $\rho$ was Corollary 6.\hfill\hal
\bigskip

Finally, we show that if $A$ is a commutative von Neumann algebra then
abelian operator modules and Finsler modules coincide.
We say that a Banach
$L^\infty(X)$-module $E$ has the {\bf $L^\infty$ norm property} if
$$
\|x\| = {\rm max}(\|p x\|, \|(1 - p) x\|)
$$
for any $x \in E$ and any projection $p \in L^\infty(X)$.
\bigskip

\noindent {\bf Theorem 9.} {\it Let $A = L^\infty(X)$ be a commutative
von Neumann algebra and let $E$ be a Banach $A$-module.
The following are equivalent:
\medskip

{\narrower{
\noindent (a) $E$ is an abelian operator $A$-module.
\medskip

\noindent (b)
There exists a map $\rho: E \to A_+$ (necessarily unique) which induces
the given norm on $E$ and makes $E$ into a Finsler $A$-module.
\medskip

\noindent (c) $E$ satisfies the $L^\infty$ norm property.
\medskip}}}

\noindent {\it Proof.}
(a) $\Rightarrow$ (b). Let $A$ and $E$ be embedded in a commutative
C*-algebra $B$; without loss of generality suppose $B$ has a unit and
$A$ is embedded unitally. (Otherwise replace $B$ by $1_A B$.
This does not alter the embedding of $E$ in $B$ since $x = 1_A x$ for
all $x \in E$.) For each $x \in E$ define
$$
\rho(x) = \inf\{f \in A: |x| \leq f\}
= \inf\{f \in A: |x| \leq f {\hbox{ and }} \|f\| = \|x\|\}.
$$
The two infima are equal since everything in the first set dominates
something in the second set.
Since $A$ is unitally embedded in $B$, the second set
contains $f = \|x\| \cdot 1$, so it is nonempty, bounded, and
self-adjoint; therefore its infimum exists.
This also shows that $\|\rho(x)\| \leq \|x\|$, and conversely, as $|x|
\leq f$ implies
$\|x\| \leq \|f\|$, we see that $\|\rho(x)\| = \|x\|$. So $\rho$ does
induce the original norm on $E$.
This automatically implies Finsler condition
(1), and condition (2) in the form $\rho(ax) = |a| \rho(x)$ is easy:
$$
\rho(ax) = \inf\{f \in A: |ax| \leq f\}
   = |a| \inf\{f \in A: |x| \leq f\} = |a| \rho(x).
$$
Again, uniqueness of $\rho$ follows from Corollary 6.
\medskip

(b) $\Rightarrow$ (c). Let $p$ be a projection in $L^\infty(X)$, let
$q = 1 - p$, and let $x \in E$. Then
$$
\rho(x) = \rho((p + q) x) = (p + q) \rho(x) = p \rho(x) + q \rho(x)
= \rho(p x) + \rho(q x).
$$
Since $\rho(px) = p \rho(x)$ and $\rho(qx) = q \rho(x)$ have disjoint
support, we get
$$
\|x\|_E = \|\rho(x)\| = {\rm max}(\|\rho(px)\|, \|\rho(qx)\|)
= {\rm max}(\|px\|_E, \|qx\|_E)
$$
as desired.
\medskip

(c) $\Rightarrow$ (a).
We assume the $L^\infty$ norm property and prove that
$E$ is $A$-convex; this suffices by Theorem 7.
\medskip

Thus let $x, y \in E$ and let $f, g \in L^\infty(X)$ be positive
functions
such that $f + g = 1$. Let $\vepsilon > 0$.
Partition $X$ into measurable subsets $X_1, \ldots, X_n$ such that $f$
and $g$ each vary by less than $\vepsilon' = \vepsilon/(\|x\| + \|y\|)$
on each $X_j$. Let $p_j$ be the characteristic function of $X_j$.
\medskip

Fix $j$ and let $\a, \b \in \R^+$ satisfy $\a + \b = 1$ and
$$
\|p_j (f - \a 1)\|, \|p_j (g - \b 1)\| \leq \vepsilon'.
$$
Then
$$
\eqalign{\|p_j (fx + gy)\| &\leq \|p_j ((f - \a 1)x + (g - \b 1)y)\|
+ \|p_j (\a x + \b y)\|\cr
&\leq \vepsilon' \|x\| + \vepsilon' \|y\| + \a \|x\| + \b \|y\|\cr
&\leq \vepsilon + {\rm max}(\|x\|, \|y\|).\cr}
$$
But then the $L^\infty$ norm property implies that
$$
\eqalign{\|fx + gy\|
&= {\rm max}(\|p_1 (fx + gy)\|, \ldots, \|p_n (fx + gy)\|)\cr
&\leq \vepsilon + {\rm max}(\|x\|, \|y\|),\cr}
$$
which in the limit $\vepsilon \to 0$ establishes that $E$ is
$A$-convex.\hfill\hal
\bigskip

A simple example of an abelian operator module which is not a Finsler
module in the C* case is given by
$A = C ([0, 1])$ and $E = L^\infty ([0, 1]).$
\bigskip
\bigskip

\noindent {\underbar{\bf 3. The noncommutative case.}}
\bigskip

\noindent {\bf Lemma 10.}
{\it Let $A$ be a C*-algebra. Then $A$ has a unique maximal commutative
ideal $I$ and it may be obtained as the intersection of the kernels of
all irreducible representations of $A$ of dimension greater than $1$.
Moreover, $I$ is contained in the center of $A$.}
\medskip

\noindent {\it Proof.} Let $I$ be the intersection of the kernels of all
irreducible representations of $A$ of dimension greater than $1$.
If $J$ is any ideal of $A$ then any irreducible representation of $A$
either annihilates $J$ or restricts to an irreducible representation of
$J$ ([4], Theorem 1.3.4). As no commutative C*-algebra has irreducible
representations of dimension greater than $1$,
it follows that $I$ contains every commutative ideal.
\medskip

Now let $a \in I$ and $b \in A$.
Then $a b - b a$ is annihilated by
every homomorphism into $\C$, and it also belongs to $I$, hence it is
annihilated by every irreducible representation of $A$.
Thus $a b - b a = 0$,
and we conclude that $I \subset Z(A).$ In particular, $I$ is
commutative.\hfill\hal
\bigskip

Recall that if $A,$ $B,$ and $D$ are C*-algebras, and if homomorphisms
$\ph : A \to D$  and $\ps : B \to D$ are given, then the C*-algebra
$A \oplus_D B$ is defined as
$$
A \oplus_D B = \{ (a, b) \in A \oplus B : \ph (a) = \ps (b) \}.
$$
We use the same notation for modules, Banach spaces, etc.

\bigskip

\noindent {\bf Lemma 11.} {\it Let $A$ and $I$ be as in Lemma 10. Then
every multiplicative linear functional on $I$ extends uniquely to a
multiplicative linear functional on $A$.
The intersection of the kernels of these functionals is an ideal $J$ of
$A$ with the properties that $I \cap J = 0$ and $A/J$ is commutative.
We have $A \cong C_0 (X) \oplus_{C_0 (Y)} B$
where $X = \Prim (A/J)$, $B = A/I$, and $Y = \Prim (A/(I + J))$.}
\medskip

\noindent {\it Proof.}
Let $\widehat I$ denote the spectrum of $I$. Every
$\omega \in \widehat I$ extends uniquely to a multiplicative linear
functional $\ep_{\omega}$ on $A$ by Theorem 1.3.4 of [4]. Let
$J = \bigcap_{\omega \in \widehat I}
\ke (\ep_{\omega})$.
\medskip

Since the range of each $\ep_{\omega}$ is $\C$, it follows that each
$\ke (\ep_{\omega})$ contains the commutator ideal of $A.$
Hence so does $J$,
so that $A/J$ is commutative.
For every nonzero element of $I$ there exists an $\omega \in \widehat I$
which does not annihilate it, so that $I \cap J = 0$.
\medskip

Define $\ph : A \to C_0 (X) \oplus_{C_0 (Y)} B$ by
$\vphi (a) = ( a + J, a + I).$ Then $\ph$ is
clearly a $*$-homomorphism, and
it is injective because $I \cap J = 0$.
To see that $\vphi$ is surjective, let $b, c \in A$ satisfy
$b + (I + J) = c + (I + J)$. Then $c - b \in I + J$ and so by the
natural
isomorphism $J = J/(I \cap J)
\cong (I + J)/I$ there exists $b' \in J$ such that $b' + I = (c - b)
+ I$.
Then $(b + b') + J = b + J$ and $(b + b') + I = c + I$,
that is, $\vphi(b + b') = (b + J, c + I)$.\hfill\hal
\bigskip

Note that $J$ can also be described in the following way.
If $U$ is the open
subset of $\Prim (A)$ corresponding to $I$, then $J$ corresponds
to the interior of the complement of $U$.
\medskip

For subsets $B$ of a C*-algebra $A$ and $F$ of a Banach $A$-module $E,$
we denote by $BF$ the closed linear span of all products $ax$ with
$a \in B$ and $x \in F.$
\bigskip

\noindent {\bf Lemma 12.} {\it Let $E$ be a Finsler module over a
C*-algebra $A$, let $I$ be an ideal of $A$, let $B = A/I$, let
$\pi: A \to 
B$ be the quotient map, and let $\rho' = \pi \circ \rho$. Then $IE =
\ke (\rho')$, $E/IE$ is a $B$-module, and $\rho'$ descends to a
$B$-valued Finsler norm on $E/IE$.}
\medskip

\noindent {\it Proof.}
It is clear that $E/IE$ is naturally
a $B$-module. If $a \in I$ and $x \in E$ then $\rho(ax)^2 =
a \rho(x)^2 a^* \in I$, hence $\rho(ax) \in I$,
that is, $ax \in \ke (\rho')$.
This shows that $IE \subset \ke (\rho')$.
Conversely, if $x \in \ke (\rho')$ then
$\rho(x) \in I$, and so there exists a sequence $\{e_n\}$ of positive
elements of $I$ such that $e_n \rho(x) \to \rho(x)$. We claim that
$$
\rho(x - e_n x)^2 = \rho(x)^2 - \rho(x)^2 e_n - e_n \rho(x)^2 +
e_n \rho(x)^2 e_n.
$$
To see this, let $b$ be the left side and $c$ the right side.
Then for any $a \in A$ we have
$$
a c a^* = (a - a e_n) \rho(x)^2 (a^* - e_n a^*) =
    \rho ((a - a e_n) x)^2 = a \rho (x - e_n x)^2 a^* = a b a^*.
$$
Thus, $a (b - c) a^* = 0$ for all $a \in A,$ whence $b = c$ as claimed.
It now follows that
$\rho(x - e_n x)^2 \to 0$. Thus $\|x - e_n x\|_E \to 0$ and so
$x \in IE$.
We have therefore shown that $IE = \ke (\rho')$.
\medskip

Next we show that $\|\rho'(\cdot)\|$ satisfies the triangle inequality.
For suppose this fails and
$$
\|\rho'(x + y)\| > \|\rho'(x)\| + \|\rho'(y)\|
$$
for some $x, y \in E$.
Then there is a pure state $f$ on $B$ such that
$f(\rho'(x + y)^2) = \alpha^2$, $f(\rho'(x)^2) = \beta^2$, and
$f(\rho'(y)^2) = \gamma^2$ with
$\alpha = \|\rho'(x + y)\|$, $\beta \leq \|\rho'(x)\|$, and $\gamma \leq
\|\rho'(y)\|$ (hence $\alpha > \beta + \gamma$).
Letting $f' = f \circ \pi$, we
get that $f'$ is a pure state on $A$ and $f'(\rho(x + y)^2) = \alpha^2$,
$f'(\rho(x)^2) = \beta^2$, $f'(\rho(y)^2) = \gamma^2$.
\medskip

By Proposition 2.2 of [3],
there exists a net $\{a_\ld\}$ of positive
norm one elements of $A$ such that
$$
\lim \|a_\ld b a_\ld - f'(b) a_\ld^2\| = 0
$$
for any $b \in A$.
It follows that $\|a_\ld b a_\ld\| \to f'(b)$ for
all $b$. In particular with $b = \rho(x + y)^2$ we have
$$
\|\rho(a_\ld x + a_\ld y)^2\| =
\|a_\ld \rho(x + y)^2 a_\ld\| \to \alpha^2,
$$
hence
$$
\|a_\ld x + a_\ld y\|_E = \|\rho(a_\ld x + a_\ld y)\|
\to \alpha,
$$
and similarly $\|a_\ld x\|_E \to \beta$ and $\|a_\ld y\|_E \to
\gamma$. Since $\alpha > \beta + \gamma$ this contradicts the triangle
inequality in $E$.
We conclude that $\|\rho'(\cdot)\|$ satisfies the triangle
inequality.
\medskip

We must now show that $\rho'$ descends to $E/IE$, that is, we must prove
that $\rho'(x) = \rho'(x + y)$ for any $x \in E$ and $y \in IE$.
The triangle inequality just established implies that for any
$a \in A_+$ we have
$$
\|\rho'(ax)\| = \|\rho'(ax)\| - \|\rho'(ay)\| \leq \|\rho'(ax + ay)\|
\leq \|\rho'(ax)\| + \|\rho'(ay)\| = \|\rho'(ax)\|,
$$
whence $\|\rho'(ax)\| = \|\rho'(ax + ay)\|.$
If $b \in B_+,$ then there is $a \in A_+$ such that $\pi (a) = b,$
so that
$$
\| b \rho'(x)^2 b \| = \|\rho'(ax) \|^2 = \|\rho'(ax + ay)\|^2
          = \| b \rho'(x + y)^2 b\|.
$$
Theorem 4 now implies that $\rho'(x) = \rho'(x + y),$ showing that
$\rho'$ does descend to $E/IE$.\hfill\hal
\bigskip

\noindent {\bf Lemma 13.} {\it Retain the notation of Lemma 12.
Suppose $\rho'$ satisfies the parallelogram law
$$
\rho'(x + y)^2 + \rho'(x - y)^2 = 2 \rho'(x)^2 + 2 \rho'(y)^2
$$
for $x, \, y \in E.$
Then $E/IE$ is a Hilbert $B$-module for a unique $B$-valued
inner product which gives rise to $\rho'$.}
\medskip

\noindent {\it Proof.} We begin by showing that the polarization formula
$$
\la x, y\ra = \f \sum_{k=0}^3 i^k \rho'(x + i^k y)^2
$$
defines a $\C$-sesquilinear map $\la \cdot,\cdot \ra: E \times E \to B$
which satisfies $\la x, y \ra^* = \la y, x \ra$ and
$\rho'(x)^2 = \la x, x \ra$ for all $x, y \in E$.
\medskip

First, we have
$$
\la x, x \ra = \f \sum_{k=0}^3 i^k \rho'((1 + i^k)x)^2
 = \f \sum_{k=0}^3 i^k (1 + i^k) \rho'(x)^2 (1 + i^{-k})
 = \rho'(x)^2.
$$
\medskip

Next, observe that $\rho'(i^k z)^2 = i^k \rho'(z)^2 i^{-k}
= \rho'(z)^2$, whence
$$
\eqalign{\la y, x \ra &= \f \sum_{k=0}^3 i^k \rho'(y + i^k x)^2
 = \f \sum_{k=0}^3 i^k \rho'(i^k (x + i^{-k} y))^2  \cr
&= \f \sum_{k=0}^3 i^{-k} \rho'(x + i^k y)^2
 = \la x, y \ra^*.}
$$
This proves one of the claims and also shows that to prove
$\C$-sesquilinearity we need only check $\C$-linearity
in the first variable.
This is done by exactly the same argument that one uses in the scalar
case (e.g.\ see [6]), but we include this argument for completeness.
\medskip

For $u,v,z \in E$ the parallelogram law gives
$$
\rho'((u + i^k z) + v)^2 + \rho'((u + i^k z) - v)^2
= 2 \rho'(u + i^k z)^2 + 2 \rho'(v)^2,
$$
that is,
$$
\rho'((u + v) + i^k z)^2 + \rho'((u - v) + i^k z)^2
= 2 \rho'(u + i^k z)^2 + 2 \rho'(v)^2.
$$
Multiplying by $i^k / 4$ and summing over $k$ yields
$$
\la u + v, z \ra + \la u - v, z \ra = 2\la u , z\ra.\eqno(*)
$$
Substituting $v = u$ then shows that $\la 2 u, z \ra = 2 \la u, z \ra$,
using
the fact that
$$
\la 0, z \ra = \f \sum_{k=0}^3 i^k \rho'(i^k z)^2 =
\f \sum_{k=0}^3 i^k \rho'(z)^2 = 0.
$$
So replacing $2 \la u, z \ra$ with $\la 2u, z \ra$ in ($*$) and
substituting
$u = (x + y)/2$ and $v = (x - y)/2$, we get
$$
\la x, z \ra + \la y, z \ra = \la x + y, z \ra.
$$
This proves additive linearity.
\medskip

Now for $\a \in \R$ define $f(\a) = \la \a x, y \ra \in A$. The map
$\a \mapsto \a x + i^k y$ is continuous since $E$ is a Banach module,
so Corollary 5 implies that the map $\a \mapsto \rho'(\a x + i^k y)^2$
is continuous.
It follows that $f$ is continuous. As we also have $f(\a +
\b) = f(\a) + f(\b)$, it follows that $f(\a) = \a f(1)$ for
all $\a \in \R$, that is, $\la \a x, y \ra = \a \la x, y \ra$.
Finally, direct calculation shows that
$$
\eqalign{\la i x, y \ra &= \f \sum_{k=0}^3 i^k \rho'(i x + i^k y)^2
 = \f \sum_{k=0}^3 i^k \rho'(x + i^{k-1} y)^2\cr
&= \f \sum_{k=0}^3 i^{k+1} \rho'(x + i^k y)^2
 = i \la x, y \ra.\cr}
$$
So additive linearity implies that
$$
\la (\a + i \b)x, y\ra = \la \a x, y\ra + \la i \b x,y\ra
= \a \la x, y \ra + i \b \la x, y \ra = (\a + i \b) \la x, y \ra.
$$
This completes the proof of $\C$-linearity in the first variable, hence
of $\C$-sesquilinearity by the comment made earlier.
\medskip

Let $y \in IE$, so that $\la y, y \ra = \rho'(y)^2 = 0$.
For any positive linear functional $f$ on $B$,
$f(\la \cdot, \cdot \ra)$ is a $\C$-valued positive semidefinite
sesquilinear form, hence it satisfies the Cauchy-Schwartz inequality. In
particular
$$
|f(\la x, y \ra)|^2 \leq f(\la x, x \ra) f(\la y, y \ra) = 0.
$$
It follows that $\la x, y \ra = 0$ for all $x \in E$, and from this we
conclude that $\la \cdot, \cdot \ra$ descends to $E/IE$.
\medskip

We now know that $\la x + IE, y + IE \ra = \la x, y \ra$
defines a $B$-valued
inner product on $E/IE$ that satisfies $\la x + IE, y + IE \ra^* =
\la y + IE, x + IE \ra$ and $\rho'(x)^2 = \la x + IE, x + IE \ra$ and is
$\C$-sesquilinear. To complete the proof that $E/IE$ is
a Hilbert $B$-module we must prove $B$-sesquilinearity. To do this
fix $x,y \in E$ and consider the $\C$-sesquilinear forms
$\{ \cdot, \cdot \}, \{ \cdot, \cdot \}': A\times A
\to B$ defined by $\{a, b\} = \pi(a) \la x, y \ra \pi(b)^*$ and
$\{a, b\}' = \la ax, by \ra$. For any $a \in A$ we have
$$
\eqalign{\{ a, a \}
&= \f \sum_{k=0}^3 i^k \pi(a) \rho'(x + i^k y)^2 \pi(a)^*
 = \f \sum_{k=0}^3 i^k \pi(a \rho(x + i^k y)^2a^*)\cr
&= \f \sum_{k=0}^3 i^k \rho'(a x + i^k a y)^2
 = \{ a, a \}'.\cr}
$$
It follows that $\{ \cdot, \cdot \} = \{ \cdot, \cdot \}'$ by
polarization.
Thus
$$
\pi(a) \la x + IE, y + IE \ra \pi(b)^* =
\la \pi(a) (x + IE), \pi(b) (y + IE) \ra
$$
for all $a,b \in A$ and $x,y \in E$, and this shows that
$\la \cdot, \cdot \ra$ is $B$-sesquilinear.
\medskip

Finally, $\la \cdot, \cdot \ra$ is unique by polarization.\hfill\hal
\bigskip

The following lemma is a much simpler relative of Lemma 6.7.1 of [17].
The proof is simple enough that we give it anyway.
\bigskip

\noindent {\bf Lemma 14.} {\it Let $A \subset B(H)$ be a C*-algebra
irreducibly represented on a Hilbert space $H$ of dimension
greater than $1$. Let $\xi \in H$.
Then there exist $a, b \in A$ such that $a = a^*$, $a \xi = \xi$,
$b a = 0$, and $b b^* = a^2$.}
\medskip

\noindent {\it Proof.}
Let $\zeta \in H$ be orthogonal to $\xi$ and satisfy
$\|\zeta\| = \|\xi\|$. By the Kadison Transitivity Theorem
([17], Theorem 2.7.5) we can find self-adjoint $c,s \in A$
such that $c \xi = \xi$, $c \zeta = 0$, $s \xi = \zeta$, and
$s \zeta = \xi$. Choose continuous functions $f,g: \R \to [0,1]$ with
$f(1) = 1$, $g(0) = 1$, and $fg = 0$. Set $b = f(c) s g(c)$. Note that
$b \zeta = \xi$ and $b^* \xi = \zeta$, so $b b^* \xi = \xi$. Set $a =
(b b^*)^{1/2}$. Then $fg = 0$ implies $b^2 = 0$ implies $ba (ba)^* =
b^2 (b^*)^2 = 0$, so $ba = 0$.\hfill\hal
\bigskip

\noindent {\bf Lemma 15.}
{\it Let $E$ be a Finsler $A$-module, let $I$ be the maximal
commutative ideal of $A$, and let $B = A/I$. Then $E/IE$ is a
Hilbert $B$-module.}
\medskip

\noindent {\it Proof.} Let $\pi: A \to B$ be the natural projection
and let $\rho' = \pi \circ \rho$. By Lemma 13 it will suffice to show
that
$\rho'$ satisfies the parallelogram law. Let $x,y \in E$, let
$\vphi: A \to B(H)$ be an irreducible representation on a Hilbert space
of dimension greater than $1$, and let $\xi \in H$.
\medskip

By Lemma 14 there exist $a,b \in A$ such that the elements
$\tilde{a} = \vphi(a)$ and $\tilde{b} = \vphi(b)$ satisfy
$\tilde{a} = \tilde{a}^*$,
$\tilde{a} \xi = \xi$, $\tilde{b} \tilde{b}^* = \tilde{a}^2$, and
$\tilde{b} \tilde{a} = 0.$
Letting $\tilde{\rho} = \vphi \circ \rho$ and using the
fact proved in Lemma 12 that $\tilde{\rho}$ descends to
$E/\ke (\tilde{\rho}) = E/\ke (\vphi)E$
(and $ba, bb^* - a^2 \in \ke (\vphi)$) we have
$$
\eqalign{\la \tilde{\rho}(x \pm y)^2 \xi, \xi \ra
&= \la \tilde{a}^2 \tilde{\rho}(x \pm y)^2 \tilde{a}^2 \xi, \xi \ra
 = \la \tilde{\rho}(a^2(x \pm y))^2 \xi, \xi \ra\cr
&= \la \tilde{\rho}((a \pm b)(a x + b^* y))^2 \xi, \xi \ra
 = \la (\tilde{a} \pm \tilde{b})\tilde{\rho}(a x + b^* y)^2(\tilde{a}
\pm \tilde{b}^*)\xi, \xi \ra.\cr}
$$
Adding yields
$$
\eqalign{ \big\langle  \big[ \tilde{\rho}(x + y)^2 +
               \tilde{\rho}(x - y)^2  \big] \xi,\xi \big\rangle
&= \big\langle \big[ 2\tilde{a} \tilde{\rho}(a x + b^* y)^2 \tilde{a} +
 2\tilde{b} \tilde{\rho}(a x + b^* y)^2 \tilde{b}^* \big]
                                                      \xi, \xi \big\rangle\cr
&=  \big\langle  \big[ 2\tilde{\rho}(a(a x + b^* y))^2 +
 2\tilde{\rho}(b(a x + b^* y))^2 \big] \xi, \xi \big\rangle\cr
&=  \big\langle  \big[ 2\tilde{\rho}(a^2 x)^2 + 2\tilde{\rho}(a^2 y)^2 \big]
                                                       \xi, \xi \big\rangle\cr
&=  \big\langle \tilde{a}^2 \big[ 2\tilde{\rho}(x)^2 +
        2\tilde{\rho}(y)^2 \big] \tilde{a}^2\xi, \xi \big\rangle\cr
&=  \big\langle  \big[ 2\tilde{\rho}(x)^2 + 2\tilde{\rho}(y)^2 \big]
                                              \xi, \xi \big\rangle.\cr}
$$
Let $c = \rho(x + y)^2 + \rho(x - y)^2 - 2\rho(x)^2 - 2\rho(y)^2$, so we
have $\la c \xi, \xi\ra = 0$.
As $\xi$ was arbitrary we get $\vphi(c) = 0$,
and since this is true for all irreducible
representations $\vphi$ of dimension greater than $1$ we
conclude that $c \in I$. This implies that $\pi(c) = 0$ and so $\rho'$
satisfies the parallelogram law, as desired.\hfill\hal
\bigskip

\noindent {\bf Lemma 16.}
{\it Let $A = B_1 \oplus_{D} B_2,$ with $\ph_i : B_i \to D$ surjective.
\medskip

(1) Let $F_1$ and $F_2$ be Finsler modules over $B_1$ and $B_2,$ let
$H$ be a Finsler module over $D,$ and let $\ps_i : F_i \to H$ be
continuous linear maps inducing Finsler module isomorphisms
$\overline{\ps}_i : F_i / \ke (\ph_i) F_i \to H.$ Then
$E = F_1 \oplus_{H} F_2$ is a Finsler module over $A,$ with the
module structure $(b_1, b_2) (x_1, x_2) = (b_1 x_1, b_2 x_2)$ and
$\rh_E (x_1, x_2) = ( \rh_{F_1} (x_1), \rh_{F_2} (x_2)).$
\medskip

(2) Let $E$ be a Finsler module over $A.$
Let $\pi_i : A \to B_i$ be the projection maps, and set
$F_i = E / \ke (\pi_i) E$
and $H = E / \ke (\ph_1 \circ \pi_1) E = E / \ke (\ph_2 \circ \pi_2) E.$
Then there is a canonical isomorphism $E \cong F_1 \oplus_{H} F_2.$}
\medskip

\noindent {\it Proof.}
(1) This is a straightforward calculation, and is omitted.
\medskip

(2) We have $\ph_1 \circ \pi_1 = \ph_2 \circ \pi_2$ by the
definition of $B_1 \oplus_{D} B_2.$ Since $\ph_1$ and $\ph_2$ are
surjective, so are $\pi_1$ and $\pi_2.$ Therefore Lemma 12 shows
that $F_i$ is a Finsler module over $B_i$ and $H$ is a Finsler module
over $D.$ It is immediate that the map
$x \mapsto (x + \ke (\pi_1) E, x + \ke (\pi_2) E)$ is a
homomorphism of
$A$-modules, and easy to check that it intertwines the Finsler norms.
Since it intertwines the Finsler norms, it must be injective,
and it is surjective by the argument used at the end of the
proof of Lemma 11.\hfill\hal
\bigskip

\noindent {\bf Theorem 17.}
{\it Let $A$ be a C*-algebra and as in Lemma 11
write $A \cong C_0(X) \oplus_{C_0(Y)} B$, where $C_0(X) = A/J$,
$B = A/I$, and $C_0(Y) = A/(I + J)$.
If $E_0$ is a Hilbert module over $C_0(Y)$, $E_1$ is a Finsler module
over $C_0(X)$, $E_2$ is a Hilbert module over $B$,
and $\ps_i : E_i \to E_0$ are
continuous linear maps inducing Finsler module isomorphisms
$$
\overline{\ps}_1 : E_1/ [(I + J)/J] E_1 \to E_0 \,\,\,\,\,\, {\rm{and}}
\,\,\,\,\,\,
\overline{\ps}_2 : E_2/ [(I + J)/I] E_2 \to E_0,
$$
then $E_1 \oplus_{E_0} E_2$ is a Finsler module over $A$.
Conversely, every
Finsler module over $A$ arises in this way.}
\medskip

\noindent {\it Proof.}
If in the statement we merely require $E_0$ and $E_2$ to be Finsler
modules, then this is just the previous lemma. However, Lemma 15
implies that $E_2$ is necessarily a Hilbert module, and it follows
that $E_0 \cong E_2/ [(I + J)/I] E_2$ is also automatically
a Hilbert module.\hfill\hal
\bigskip

\noindent {\bf Corollary 18.}
{\it Let $A$ be a C*-algebra. Then the class of Finsler $A$-modules
equals the class of Hilbert $A$-modules if and only if $A$ has no
nonzero commutative ideals. In particular this holds if $A$ is simple
with $\dim (A) > 1,$
approximately divisible [7], or a von Neumann algebra with no abelian
summand.}
\medskip

\noindent {\it Proof.} We use the notation of Theorem 17.
If $A$ has no commutative ideals then $I = 0$,
$J = A$, $C_0(X) = C_0(Y) = 0$, and $B = A$.
Thus Theorem 17 identifies the Finsler modules over $A$ with the Hilbert
modules over $B = A$. Conversely, if $I$ is a nontrivial commutative
ideal of $A$ then $C_0(X) \neq 0.$
Choose a non-Hilbert Banach space $V,$ and
take $E_1$ to be the module of continuous maps
$\widehat I \to V$ which
vanish at infinity. Letting
$E_2 = E_0 = 0$, Theorem 17 produces a Finsler $A$-module; but it is not
a Hilbert $A$-module because the norm $\rho$ does not satisfy the
parallelogram law.\hfill\hal

\bigskip
\bigskip

[1] C.\ A.\ Akemann, personal communication.
\medskip

[2] C.\ A.\ Akemann, J.\ Anderson, and G.\ K.\ Pedersen, Triangle
inequalities in operator algebras, {\it Lin.\ and Mult.\ Algebra \bf 11}
(1982), 167-178.
\medskip

[3] ----------, Excising states of C*-algebras, {\it Can.\ J.\ Math.\
\bf 38} (1986), 1239-1260.
\medskip

[4] W.\ B.\ Arveson, {\it An Invitation to C*-algebras}, Springer GTM
{\bf 39} (1976).
\medskip

[5] D.\ Bao and S.\ S.\ Chern, A note on the Gauss-Bonnet theorem for
Finsler spaces, {\it Ann.\ of Math.\ \bf 143} (1996), 233-252.
\medskip

[6] S.\ K.\ Berberian, {\it Lectures in Functional Analysis and Operator
Theory}, Springer GTM {\bf 15} (1974).
\medskip

[7] B.\  Blackadar, A.\  Kumjian, and  M.\  R\o rdam,
Approximately central matrix units and the structure of
noncommutative tori, {\it K-Theory {\bf 6}} (1992), 267-284.
\medskip

[8] D.\ P.\ Blecher, A new approach to Hilbert C*-modules,
{\it Math.\ Ann.}, to appear.
\medskip

[9] E.\ Cartan, {\it Les Espaces de Finsler},
Actualit\'es Scientifiques et Industrielles no.\ 79 (1934).
\medskip

[10] J.\ Dauns and K.\ H.\ Hoffman, {\it Representation of Rings
by Sections}, {\it Mem.\ Amer.\ Math.\ Soc.\ \bf 83} (1968).
\medskip

[11] K.\ de Leeuw, Banach spaces of Lipschitz functions,
{\it Studia Math.\  \bf 21} (1961), 55-66.
\medskip

[12] M.\ J.\ Dupr\'e and R.\ M.\ Gillette, {\it Banach bundles, Banach
modules and automorphisms of C*-algebras}, Pitman Research Notes in
Mathematics {\bf 92} (1983).
\medskip

[13] E.\ G.\ Effros and Z.-J.\ Ruan,
Representations of operator bimodules
and their applications, {\it J.\ Op.\ Thy.\ \bf 19} (1988), 137-157.
\medskip

[14] M.\ Frank, Geometrical aspects of Hilbert C*-modules, preprint.
\medskip

[15] R.\ Godement, Th\'eorie g\'en\'erale des sommes continus d'espaces de
Banach, {\it C.\ R.\ Acad.\ Sci.\ Paris \bf 228} (1949), 1321-1323.
\medskip

[16] I.\ Kaplansky, Modules over operator algebras,
{\it Amer.\ J.\ Math.\ \bf 75} (1953), 839-858.
\medskip

[17] G.\ K.\ Pedersen, {\it C*-algebras and their Automorphism Groups},
Academic Press (1979).
\medskip

[18] M.\ A.\ Rieffel, Morita equivalence for operator algebras,
{\it Proc.\ Symp.\ Pure Math.\ \bf 38} (1982), 285-298.
\medskip

[19] J.-L.\ Sauvageot, Tangent bimodule and locality for dissipative
operators on C*-algebras, {\it Quantum Probability and App.\ IV},
Springer LNM {\bf 1396} (1989), 322-338.
\medskip

[20] R.\ G.\ Swan, Vector bundles and projective modules, {\it Trans.\
Amer.\ Math.\ Soc.\ \bf 105} (1962), 264-277.
\medskip

[21] A.\ Takahashi, {\it Fields of Hilbert Modules}, Dissertation,
Tulane University (1971).
\medskip

[22] N.\ Weaver,
Lipschitz algebras and derivations of von Neumann algebras,
{\it J.\ Funct.\ Anal.\ \bf 139} (1996), 261-300.
\medskip

[23] ----------, Deformations of von Neumann algebras,
{\it J.\ Op.\ Thy.\ \bf 35} (1996), 223-239.
\medskip

[24] ----------, Operator spaces and noncommutative metrics, preprint.
\medskip

[25] ----------, $\a$-Lipschitz algebras on the noncommutative torus,
preprint.
\bigskip
\bigskip

\noindent Dept.\ of Mathematics

\noindent University of Oregon

\noindent Eugene, OR 97403

\noindent phillips@math.uoregon.edu
\bigskip
\bigskip

\noindent Dept.\ of Mathematics

\noindent UCLA

\noindent Los Angeles, CA 90024

\noindent nweaver@math.ucla.edu
\end